\documentclass{amsart}
\begin{document}
\title[$p$-jets]{The first $p$-jet space  of an elliptic curve:\\
global functions and lifts of Frobenius}
\bigskip

\def \Rp{R_p}
\def \Rpi{R_{\pi}}
\def \dpi{\d_{\pi}}
\def \bT{{\bf T}}
\def \cI{{\mathcal I}}
\def \cH{{\mathcal H}}
\def \cJ{{\mathcal J}}
\def \ZN{\bZ[1/N,\zeta_N]}
\def \tA{\tilde{A}}
\def \o{\omega}
\def \tB{\tilde{B}}
\def \tC{\tilde{C}}
\def \alph{A}
\def \bet{B}
\def \bsigma{\bar{\sigma}}
\def \y{^{\infty}}
\def \Ra{\Rightarrow}
\def \uBS{\overline{BS}}
\def \lBS{\underline{BS}}
\def \lB{\underline{B}}
\def \<{\langle}
\def \>{\rangle}
\def \hL{\hat{L}}
\def \cU{\mathcal U}
\def \cF{\mathcal F}
\def \S{\Sigma}
\def \st{\stackrel}
\def \sd{Spec_{\d}\ }
\def \pd{Proj_{\d}\ }
\def \s{\sigma_2}
\def \i{\sigma_1}
\def \bs{\bigskip}
\def \cD{\mathcal D}
\def \cC{\mathcal C}
\def \cT{\mathcal T}
\def \cK{\mathcal K}
\def \cX{\mathcal X}
\def \sX{X_{set}}
\def \cY{\mathcal Y}
\def \cS{X}
\def \cR{\mathcal R}
\def \cE{\mathcal E}
\def \tcE{\tilde{\mathcal E}}
\def \cP{\mathcal P}
\def \cA{\mathcal A}
\def \cV{\mathcal V}
\def \cM{\mathcal M}
\def \cL{\mathcal L}
\def \cN{\mathcal N}
\def \tcM{\tilde{\mathcal M}}
\def \caS{\mathcal S}
\def \cG{\mathcal G}
\def \cB{\mathcal B}
\def \tG{\tilde{G}}
\def \cF{\mathcal F}
\def \h{\hat{\ }}
\def \hp{\hat{\ }}
\def \tS{\tilde{S}}
\def \tP{\tilde{P}}
\def \tA{\tilde{A}}
\def \tX{\tilde{X}}
\def \tcS{\tilde{X}}
\def \tT{\tilde{T}}
\def \tE{\tilde{E}}
\def \tV{\tilde{V}}
\def \tC{\tilde{C}}
\def \tI{\tilde{I}}
\def \tU{\tilde{U}}
\def \tG{\tilde{G}}
\def \tu{\tilde{u}}
\def \chu{\check{u}}
\def \tx{\tilde{x}}
\def \tL{\tilde{L}}
\def \tY{\tilde{Y}}
\def \d{\delta}
\def \e{\chi}
\def \bW{\mathbb W}
\def \bV{{\mathbb V}}
\def \bF{{\bf F}}
\def \bE{{\bf E}}
\def \bC{{\bf C}}
\def \bO{{\bf O}}
\def \bR{{\bf R}}
\def \bA{{\bf A}}
\def \bB{{\bf B}}
\def \cO{\mathcal O}
\def \ra{\rightarrow}
\def \bx{{\bf x}}
\def \f{{\bf f}}
\def \bX{{\bf X}}
\def \bH{{\bf H}}
\def \bS{{\bf S}}
\def \bF{{\bf F}}
\def \bN{{\bf N}}
\def \bK{{\bf K}}
\def \bE{{\bf E}}
\def \bB{{\bf B}}
\def \bQ{{\bf Q}}
\def \bd{{\bf d}}
\def \bY{{\bf Y}}
\def \bU{{\bf U}}
\def \bL{{\bf L}}
\def \bQ{{\bf Q}}
\def \bP{{\bf P}}
\def \bR{{\bf R}}
\def \bC{{\bf C}}
\def \bD{{\bf D}}
\def \bM{{\bf M}}
\def \bZ{{\mathbb Z}}
\def \xtoleqr{x^{(\leq r)}}
\def \hU{\hat{U}}
\def \k{\kappa}
\def \ee{\overline{p^{\k}}}

\newtheorem{THM}{{\!}}[section]
\newtheorem{THMX}{{\!}}
\renewcommand{\theTHMX}{}
\newtheorem{theorem}{Theorem}[section]
\newtheorem{corollary}[theorem]{Corollary}
\newtheorem{lemma}[theorem]{Lemma}
\newtheorem{proposition}[theorem]{Proposition}
\theoremstyle{definition}
\newtheorem{definition}[theorem]{Definition}
\theoremstyle{remark}
\newtheorem{remark}[theorem]{Remark}
\newtheorem{example}[theorem]{\bf Example}
\numberwithin{equation}{section}
\maketitle

\bigskip

\medskip
\centerline{\bf Alexandru Buium${}^1$, Arnab Saha${}^2$}

\medskip

\noindent {\it ${}^1$University of New Mexico, Albuquerque,  NM 87131, USA; buium@math.unm.edu}
\smallskip
{\it ${}^2$Australian National University, ACT 2601, Australia; arnab.saha@anu.edu.au}

\bigskip

\begin{abstract}
We prove that there are no non-constant global functions and no lifts of Frobenius  on the first $p$-jet space on an elliptic curve unless the elliptic curve itself has a lift of Frobenius.
\end{abstract}

\section{Introduction and main results}

The aim of this paper is to settle two basic issues in the theory of arithmetic differential equations that were (implicitly) left open in the papers on the subject, in particular in \cite{char, book, BYM}. Our main result is Theorem \ref{main} below. To explain our result let us recall some basic concepts from loc. cit. Let $p$ be an odd prime and let $R$ be the $p$-adic completion of the maximum unramified extension of the ring $\bZ_p$ of $p$-adic integers. 
By a $p$-formal scheme (over $R$) we understand a formal scheme over $R$ with ideal of definition generated by $p$. A $p$-formal scheme over $R$ is said to have  a lift of Frobenius if it possesses an endomorphism  over $\bZ_p$ whose reduction mod $p$ is the absolute $p$-power Frobenius. 
For any $R$-scheme of finite type $X$ we defined in \cite{char} a projective system of $p$-formal schemes $J^n(X)$, $n\geq 0$, called the $p$-jet spaces of $X$; morally their rings of global functions $\cO(J^n(X))$ should be viewed as rings of {\it arithmetic differential equations} on $X$ \cite{book}.  If $X$ is a group scheme then $J^n(X)$ are groups in the category of $p$-formal schemes.
Cf. section 2 below for a review of these spaces. 
Here is our main result; in its statement, by an elliptic curve over $R$ we understand a smooth projective curve of genus $1$ with an $R$-point.

\begin{theorem}
\label{main}
Let $E$ be an elliptic curve over $R$ that does not have a lift of Frobenius. The following hold:

1) $J^1(E)$ has no global functions except the elements of $R$; i.e. $\cO(J^1(E))=R$;

2) $J^1(E)$ has no lift of Frobenius.
\end{theorem}

\begin{remark}
The space $J^1(E)$ played a key role in the theory of arithmetic differential equations \cite{char, book, BYM}; to give just one recent example, this space was used in \cite{BYM} as the natural phase space for the arithmetic Painev\'{e} VI equation in Hamiltonian description. Somehow the above mentioned theory could be developed by finding ways around proving a result such as Theorem \ref{main}; this theorem, however, substantially  clarifies the situation. 
\end{remark}

\begin{remark}
To place the above theorem in context let us mention that if $X$ is a smooth projective curve of genus $g$ over $R$ then the following results were previously known. 

1) If $g=0$ then $\cO(J^n(X))=R$ for $n\geq 0$; cf. \cite{pjets}. 

2) If $g\geq 2$ then 
$J^n(X)$ is affine for $n\geq 1$ \cite{pjets}; in particular $\cO(J^1(X))\neq R$ and $J^1(X)$ has a lift of Frobenius. 

3) If $g=1$
(so $X=E$, an elliptic curve) the situation is as follows:

 If $E$ {\it does} have a lift of Frobenius then, trivially, $J^1(E)=\widehat{E}\times \widehat{{\mathbb G}_a}$  as groups (where ${\mathbb G}_a$ is the additive group) and hence 
 $\cO(J^1(E))=\cO(\widehat{{\mathbb G}_a})\neq R$; also
   $J^1(E)$ possesses lifts of Frobenius (because both $E$ and ${\mathbb G}_a$ do).
   
 If  $E$ does {\it not} have a lift of Frobenius there is no non-zero homomorphism $J^1(E)\ra \widehat{{\mathbb G}_a}$ (no {\it $\d$-character of order $1$}); cf. \cite{book}, p. 194, Proposition 7.15. On the other hand there exists a non-zero homomorphism
   $\psi=\psi_{2,can}:J^2(E)\ra \widehat{{\mathbb G}_a}$ with the property that any other such homomorphism is an $R$-multiple of $\psi$. Cf.
   \cite{char} and \cite{book} (p. 201, Definition 7.24 and p. 205, Theorem 7.34). This $\psi$ can be called a {\it canonical $\d$-character of order $2$ of $E$} and should be viewed as an arithmetic analogue of the ``Manin map" of elliptic curves defined over a differential field
   \cite{manin, hermann}.  (This $\psi$ is then unique up to multiplication by an element in $R^{\times}$; given an invertible $1$-form on $E$, $\psi$ can be further normalized and is unique with this normalization; we will not need this normalized version here but this is sometimes important, e.g. in \cite{BYM}.) Although Theorem \ref{main} does not refer to $J^2(E)$ its proof will use $J^2(E)$ and the map $\psi$.
   
   4) If $E$ is an elliptic curve over $R$ we have an extension  $$0\ra \widehat{{\mathbb G}_a}\ra J^1(E) \ra \widehat{E}\ra 0$$
     of groups in the category of  $p$-formal schemes; cf. \cite{char}. Note that
   $J^1(E)$ is not, in general,  the universal vectorial extension  of $E$ in the sense of \cite{mm} i.e. the map
   $Hom(\widehat{{\mathbb G}_a},\widehat{{\mathbb G}_a})\ra H^1(E,\cO)$ is not  generally an isomorphism. However, for $E$ ``sufficiently generic over $R$" the latter map {\it is} an isomorphism. The genericity required for that is that the value of the differential modular form $f^1$ in \cite{book}, pp. 193-194,  be invertible in $R$.
   In this latter case our Theorem \ref{main} becomes, therefore, a theorem about the universal vectorial extension of $E$.
\end{remark}

\begin{remark}
Assertion 1) in Theorem \ref{main} is analogous to the situation in differential algebra  where the differential algebraic first jet space (still denoted in this remark by) $J^1(E)$ of an elliptic curve $E$ defined over a differential field $F$ (and not defined over the constants of that field) has no other global functions except the elements of $F$; this follows from the fact that $J^1(E)$ is isomorphic, in this case,  to the universal vectorial extension of $E$. Cf. \cite{hermann}, p. 125. On the other hand assertion 2) in Theorem \ref{main} is in deep contrast with the situation in differential algebra where $J^1(E)$  always possesses a derivation on its structure sheaf that lifts the derivation on $F$; this is 
because, in this case, $J^2(E)\ra J^1(E)$, has a section; the latter, in its turn, follows again trivially from the fact that $J^1(E)$ is the universal vectorial extension of $E$.
\end{remark}

{\bf Acknowledgement}. The work on this paper was partially done while the first author was visiting the Australian National University in Canberra; special thanks go to James Borger for his hospitality and for inspiring conversations. 

\section{Review of $p$-jets and complements}

We start by recalling some concepts from \cite{char, pjets, book}.
If $A$ is a ring and $B$ is an $A$-algebra then a $p$-{\it derivation} from $A$ to $B$ is a map of sets $\d:A\ra B$ such that the map $A\ra W_2(B)$, $a\mapsto (a,\d a)$ is a ring homomorphism, where $W_2(B)$ is the ring of Witt vectors on $B$ of length $2$.
For any such $\d$ the map $\phi:A\ra B$, $\phi(a)=a^p\cdot 1+p\d a$ is a ring homomorphism. For any ring $A$ we denote by $\widehat{A}$ the $p$-adic completion of $A$; a ring is called complete if it is its own completion. For any scheme $X$ we denote by $\widehat{X}$ the $p$-adic completion of $X$.
If $R=\widehat{\bZ_p^{ur}}$ is as in the introduction then $R$ has a unique $p$-derivation $\d$. We set $k=R/pR$; it is an algebraic closure of the prime field with $p$ elements.
For $A$ an $R$-algebra  and $X$ a scheme or a $p$-formal scheme over $R$,  we write $\overline{A}=A/pA=A\otimes_R k$ and $\overline{X}=X\otimes_R k$ respectively. 
All $p$-derivations on $R$-algebras are assumed compatible with that of $R$.
A {\it prolongation sequence} of rings 
(over $R$) is an inductive system of $R$-algebras $(A^n)_{n\geq 0}$
 together with $p$-derivations $\d_n:A^{n}\ra A^{n+1}$ such that $\d_n$ restricted to $A^{n-1}$ is $\d_{n-1}$. We usually write $\d_n=\d$ for all $n$. Recall the definition of $p$-formal scheme from the introduction. A prolongation sequence of $p$-formal schemes is a projective system of $p$-formal schemes $(X^n)_{n\geq 0}$ such that for any open set $U\subset X^0$ the inductive system $(\cO(U^n))$ (where $U^n$ is the pre-image of $U$ in $X^n$) has a structure of prolongation sequence of rings, and these structure are compatible when $U$ varies. Prolongation sequences form a category in the obvious way. For any scheme of finite type $X$ (or any $p$-formal scheme locally isomorphic to the $p$-adic completion of such a scheme) there exists a (necessarily unique up to isomorphism) prolongation sequence $(J^n(X))_{n\geq 0}$ with $J^0(X)=\widehat{X}$ and with the property that for any prolongation sequence $(Y^n)$ and any morphism $Y^0\ra J^0(X)$ there exists a unique morphism of prolongation sequences from $(Y^n)$ to $(J^n(X))$ 
 extending $Y^0\ra J^0(X)$. The $p$-formal schemes $J^n(X)$ are called the $p$-jet spaces of $X$. If $X$ is smooth over $R$ so are $J^n(X)$; here, for $p$-formal schemes, smoothness means ``locally isomorphic to the $p$-adic completion of a smooth scheme over $R$". 
 Actually more is true: if $X$ is affine and has an \'{e}tale map $T:X\ra {\mathbb A}^N$ then, viewing $T$ as an $N$-tuple of elements in $\cO(X)$, we have that $\cO(J^n(X))\simeq
 \cO(X)[T',...,T^{(n)}]\h$ where $T',...,T^{(n)}$ are $N$-tuples of variables that are mapped to
 $\d T,...,\d^n T$ respectively.
 For all of the above the references are \cite{char, pjets, book}. In what follows we add some complements to the above.

Let $X$ be a smooth scheme or $p$-formal scheme  over $R$. We may consider the prolongation sequence
\begin{equation}
\label{aa}
...\ra J^n(X)\ra J^{n-1}(X)\ra ... \ra J^2(X)\ra J^1(X)
\end{equation}
and also the prolongation sequence
\begin{equation}
\label{bb}
...\ra J^n(J^1(X))\ra J^{n-1}(J^1(X))\ra ... \ra J^1(J^1(X))\ra J^1(X)
\end{equation}
By the universality property of the latter there is a morphism of prolongation sequences from 
\ref{aa} to \ref{bb}; in particular there is a canonical  morphism
\begin{equation}
\label{cc}
c:J^2(X)\ra J^1(J^1(X)).
\end{equation}
Cf. \cite{BYM} where this played a key role by analogy with classical mechanics. 
We also denote by $\pi^i_X:J^i(X)\ra J^{i-1}(X)$ the canonical projections. Note that $\pi^1_X$ induces, by functoriality, a morphism $J^1(\pi^1_X):J^1(J^1(X))\ra J^1(X)$. On the other  hand
we have a morphism $\pi^1_{J^1(X)}:J^1(J^1(X))\ra J^1(X)$. We get an induced morphism
\begin{equation}
\label{vint}
\pi^{11}_X:=\pi^1_{J^1(X)}\times J^1(\pi^1_X):J^1(J^1(X))\ra J^1(X)\times_{J^0(X)} J^1(X).
\end{equation}
Assume in the following Lemma that $X$ is smooth over $R$. Also $\Delta$ below denotes the diagonal embedding.

\begin{lemma}
\label{cartesian}
The following diagram is commutative and cartesian:
$$
\begin{array}{ccc}
J^1(J^1(X)) & \stackrel{\pi^{11}_X}{\longrightarrow} & J^1(X)\times_{J^0(X)} J^1(X)\\
c \uparrow & \  & \uparrow \Delta\\
J^2(X) & \stackrel{\pi^2_X}{\longrightarrow} & J^1(X)\end{array}
$$
\end{lemma}

{\it Proof}.
It is enough to prove the above in case $X$ possesses an \'{e}tale map $T:X\ra {\mathbb A}^N$ to an affine space. View $T$ as an $N$-tuple of elements in $\cO(X)$. Denote by 
$$\begin{array}{l}
\d:\cO(X)\h\ra \cO(J^1(X))=\cO(X)[\d T]\h,\\
\ \\
$$\d_1:\cO(J^1(X))=\cO(X)[\d T]\h\ra \cO(J^1(J^1(X)))=
\cO(X)[\d T, \d_1T,\d_1\d T]\h\end{array}$$
the universal $p$-derivations; and recall also that 
$$\cO(J^2(X))=\cO(X)[\d T,\d^2T]\h;$$
cf. \cite{book}, p. 75, Proposition 3.13. Then the diagram in the statement of the lemma is deduced from the following diagram of $\cO(X)\h$-algebras:
\begin{equation}
\label{bine}
\begin{array}{ccc}
\cO(X)[\d T, \d_1T,\d_1\d T]\h & \leftarrow & \cO(X)[\d T] \widehat{\otimes}_{\cO(X)}\cO(X)[\d T]\\
\downarrow & \  & \downarrow\\
\cO(X)[\d T,\d^2T]\h & \leftarrow & \cO(X)[\d T]\h
\end{array}
\end{equation}
where the top horizontal arrow sends $\d T\otimes 1\mapsto \d T$, $1\otimes \d T\mapsto \d_1T$, the bottom horizontal arrow is the inclusion, the left vertical arrow sends $\d_1 T\mapsto \d T$, $\d_1\d T\mapsto \d^2 T$, and the right vertical arrow sends  $\d T\otimes 1\mapsto \d T$,
 $1 \otimes \d T\mapsto \d T$. The diagram \ref{bine} is clearly commutative and cartesian; this ends the proof.
\qed

\bigskip

Assume in addition that $X=G$ is a smooth group scheme. Define
\begin{equation}
L_{\d}(G)=Ker(J^1(G)\ra J^0(G)).
\end{equation}
This plays the role, in our theory, of arithmetic analogue of the Lie algebra of $G$. Cf. \cite{adel2}.
However, if $G$ is non-commutative, $L_{\d}(G)$ is also non-commutative so we denote by $1$ its identity element; if $G$ is commutative so is $L_{\d}(G)$ and its identity element will be denoted by $0$. 
The quotient map
$$J^1(G)\times J^1(G)\ra J^1(G),\ \ (a,b)\mapsto ba^{-1}$$
induces a morphism 
$$q:J^1(G)\times_{J^0(G)} J^1(G)\ra L_{\d}(G).$$
Composing the latter with the morphism
$$
\pi^{11}_G:J^1(J^1(G))\ra J^1(G)\times_{J^0(G)} J^1(G)$$
defined in \ref{vint}
we get a morphism 
\begin{equation}
\label{alpha}
l^{11}\d:J^1(J^1(G))\ra L_{\d}(G).\end{equation}
The notation $l^{11}\d$ is motivated by the analogy with Kolchin's logarithmic derivative in differential algebra \cite{kolchin}.
Clearly $\pi^{11}_G$ is a group homomorphism.
If in addition $G$ is commutative then $q$ is also a group homomorphism so
$l^{11}\d$ is a group homomorphism. Assume again $G$ not necessarily commutative.
Then note that the restriction of $l^{11}\d$ to $L_{\d}(J^1(G))$ is the map
$L_{\d}(\pi^1_G):L_{\d}(J^1(G))\ra L_{\d}(G)$ (induced by $J^1(\pi^1_G)$) composed with the antipode $L_{\d}(G)\ra L_{\d}(G)$. Moreover,
  by Lemma \ref{cartesian}, we have: 

\begin{corollary}
\label{cartesian2}
$l^{11}\d^{-1}(1)=J^2(G)$.
\end{corollary}

\begin{remark}
\label{rem}
In our proof of Theorem \ref{main} we will repeatedly use the following facts which are easily checked by passing to affine open covers. Let $Z$ be a smooth $p$-formal scheme; in particular $Z$ can be $J^n(X)$ where $X$ is a smooth $R$-scheme. Then:

1) $\cO(Z)$ are flat over $R$ and $p$-adically complete;

2) The natural map $\overline{\cO(Z)}\ra \cO(\overline{Z})$ is injective;

3) $\cO(Z\times \widehat{{\mathbb A}^1})=\cO(Z)[x]\h$, where ${\mathbb A}^1=Spec\ R[x]$.
\end{remark}

We will also need the following lemmas:

\begin{lemma}\label{l1}
Let $A\ra B$ be a homomorphism between  $p$-adically complete rings in which $p$ is a non-zero divisor. If  the induced map $\overline{A}\ra \overline{B}$ is injective (respectively finite) then the map $A\ra B$ is injective (respectively finite). 
\end{lemma}

{\it Proof}. A trivial exercise.\qed

\begin{lemma}
\label{l2}
Let $A\subset B$ be an integral extension of integral domains of characteristic zero with $A$ integrally closed. Assume there exists a derivation $\theta:B\ra B$ such that 

1) $\theta A\subset A$ and 

2) $A$ and $B$ have the same constants with respect to $\theta$; i.e. if $c\in B$ and $\theta c=0$ then $c\in A$.

Then for any  $b\in B$ with $\theta b\in A$ we have $b\in A$.
\end{lemma}

{\it Proof}.
Let $K\subset L$ be the extension of the corresponding fields of fractions and let $\theta:L\ra L$ be the unique derivation extending $\theta:B\ra B$; clearly $\theta K\subset K$. Let $f(t)=t^n+a_1t^{n-1}+...+a_n\in K[t]$ be the minimal polynomial of $b$ over $K$.
Then $\theta(f(b))=0$ hence
$$(n\theta b+\theta a_1)b^{n-1}+((n-1)a_1\theta b+\theta a_2)b^{n-2}+...=0.$$
By minimality of $f$ we get $n\theta b+\theta a_1=0$ hence 
$\theta(nb+a_1)=0$. By condition 2) we have $nb+a_1\in K$ hence
(since $K$ has characteristic zero)
$b\in K$. Since $A$ is integrally closed, $b\in A$.
\qed

\section{Proof of Theorem \ref{main}}

{\it Proof of assertion 1}. 

\medskip

{\it Step 1}. Let ${\mathbb G}_a$ be the line ${\mathbb A}^1=Spec\ R[x]$ with group structure defined by $x\mapsto x\otimes 1+1\otimes x$. Then recall from \cite{book}, p. 127, that 
\begin{equation}
\label{unu}
L_{\d}(E)\simeq \widehat{{\mathbb G}_a}=Spf\ R[x]\h;
\end{equation}
we view the above isomorphism as an identification so we have a closed immersion $\iota:\widehat{{\mathbb G}_a}\subset J^1(E)$. Let $A=\cO(J^1(E))$. We then have a restriction map
\begin{equation}
\label{rho}
\rho: A\ra R[x]\h
\end{equation}
induced by the inclusion $\iota$.
Let $\overline{G}$ be the image of the multiplication by $p$ map
$[p]:\overline{J^1(E)}\ra \overline{J^1(E)}$. Then the intersection $\overline{G}\cap \overline{{\mathbb G}_a}$ in $\overline{J^1(E)}$ is  finite over $k$. So the addition map $\beta:\overline{G}\times \overline{{\mathbb G}_a}\ra \overline{J^1(E)}$ is an isogeny and $\overline{G}\ra \overline{E}$ is an isogeny. In particular $\overline{G}$ is an elliptic curve. Using  K\"{u}nneth's formula we get a map
\begin{equation}
\label{doi}
\beta^*:\cO(\overline{J^1(E)})\ra \cO(\overline{G}\times \overline{{\mathbb G}_a})=\cO(\overline{G})\otimes \cO(\overline{{\mathbb G}_a})=k[x].
\end{equation}
One trivially checks that the map \ref{doi} coincides with the restriction map 
\begin{equation}
\label{trei}
\overline{\iota}^*:\cO(\overline{J^1(E)})\ra \cO(\overline{{\mathbb G}_a})=k[x]
\end{equation}
induced by the inclusion $\overline{\iota}:\overline{{\mathbb G}_a}\subset \overline{J^1(E)}$. We conclude that the restriction map $\overline{\iota}^*$ is injective. By Remark \ref{rem}, 2), the natural map 
$\overline{A}\ra \cO(\overline{J^1(E)})$ is injective. So the composition of the latter with $\overline{\iota}^*$
in \ref{trei} is an injective map 
\begin{equation}
\label{patru}
\overline{\rho}:\overline{A}\ra \cO(\overline{{\mathbb G}_a})=k[x].
\end{equation}
Since, by Remark \ref{rem}, 1), $A$ is $p$-adically complete and flat over $R$ it follows that $\rho$ in \ref{rho} is injective, with torsion free cokernel. 

\medskip

{\it Step 2}. If $\overline{A}=k$ then, since $A$ is flat over $R$, we get by induction that the map $R/p^nR\ra A/p^nA$ is an isomorphism for all $n$. Since $A$ is $p$-adically complete we get $A=R$ and assertion 1) in Theorem \ref{main} follows. So we may (and will) assume, in what follows, that $\overline{A}\neq k$ and we seek a contradiction. Now since $\overline{A}\neq k$ it follows that $k[x]$ is finite over $\overline{A}$. In particular $\overline{A}$ is finitely generated over $k$. Let $\overline{X}=Spec\ \overline{A}$ and let $\overline{\eta}:\overline{{\mathbb G}_a}\ra \overline{X}$ be the map induced by the map $\overline{\rho}$ in \ref{patru}. So $\overline{\eta}$ is a finite dominant map, hence it is surjective. Similarly consider the affine (a priori not necessarily Noetherian) $p$-formal scheme  $X=Spf\ A$ and the map 
\begin{equation}
\label{pii}
\eta:\widehat{{\mathbb G}_a}\ra X
\end{equation}
 induced by the map $\rho$ in \ref{rho}. 
 
 \medskip

{\it Step 3}. Consider the action by translation
\begin{equation}
\label{action}
\widehat{{\mathbb G}_a}\times J^1(E)\ra J^1(E).
\end{equation}
Using Remark \ref{rem}, 3), we get an induced map
\begin{equation}
\label{coaction}
A=\cO(J^1(E))\ra \cO(\widehat{{\mathbb G}_a}\times J^1(E))=A[x]\h.
\end{equation}
One immediately checks that \ref{coaction} induces an action 
\begin{equation}
\label{mu}
\mu:\widehat{{\mathbb G}_a}\times X\ra X;
\end{equation}
to prove coassociativity one uses, again, Remark \ref{rem}, 3).
Also clearly $\widehat{{\mathbb G}_a}$ acts on itself by translation and the map 
$
\eta:\widehat{{\mathbb G}_a}\ra X$ is equivariant because the action of $\widehat{{\mathbb G}_a}$ on itself is compatible with the action 
of $\widehat{{\mathbb G}_a}$ on $J^1(E)$ and hence with the action $\mu$ of $\widehat{{\mathbb G}_a}$ on $X$.

\medskip

{\it Step 4}. By Step 3, $\overline{{\mathbb G}_a}$ acts on $\overline{X}$ and $\overline{\eta}:\overline{{\mathbb G}_a}\ra \overline{X}$ is equivariant. Since $\overline{\eta}$ is surjective the action of $\widehat{{\mathbb G}_a}$  on $\overline{X}$ is transitive so $\overline{X}$ is a smooth affine curve. Since the group $\cO(\overline{X})^{\times}$ of invertible global functions on $\overline{X}$ injects into $\cO(\overline{{\mathbb G}_a})^{\times}=k^{\times}$ it follows that $\cO(\overline{X})^{\times}=k^{\times}$ so $\overline{X}\simeq \overline{{\mathbb A}^1}$. So $\overline{A}$ identifies via $\overline{\rho}$ with a subring $k[\overline{s}]$ of $k[x]$ where $\overline{s}=\overline{s}(x)\in k[x]$ is some polynomial. We may assume $\overline{s}(0)=0$. Let $s\in A$ be any lift of $\overline{s}$ which, viewed as an element of $R[x]\h$ satisfies $s(0)=0$. Since $A$ is flat over $R$ we get, by induction, that the natural maps
$R[s]/p^nR[s]\ra A/p^nA$ are isomorphisms for all $n$. Since $A$ is $p$-adically complete we get an isomorphism $R[s]\h\simeq A$ which we view from now on as an equality; hence $A=R[s]\h\subset R[x]\h$, so $s=s(x)$ is a restricted power series. Note that by Lemma \ref{l1} the extension $A\subset R[x]\h$ is integral. Also $A$ is integrally closed because it is a regular ring.

\medskip

{\it Step 5}. The action of $\widehat{{\mathbb G}_a}$ on $X$ is given by a map
$\mu^*:A\ra A[x]\h$,
$$\mu^* (a)=a+(Da)x+...$$
where $D:A\ra A$ is an $R$-derivation. 
By the equivariance of $\eta$ we have a commutative diagram
$$
\begin{array}{rcl}
A & \stackrel{\mu^*}{\longrightarrow} & A[x]\h\\
\rho \downarrow & \  & \downarrow \rho\otimes 1\\
R[x]\h & \longrightarrow & R[x]\h\widehat{\otimes}R[x]\h
\end{array}
$$
where $\rho$ is the inclusion and the bottom arrow sends $x\mapsto x\otimes 1+1\otimes x$, hence sends any $g\in R[x]\h$ into $g\otimes 1+ \frac{dg}{dx}\otimes x+....$.
In particular  we get
$$a\otimes 1+(Da)\otimes x+...=a\otimes 1+ \frac{da}{dx}\otimes x+...$$
for $a\in A$; so $Da=\frac{da}{dx}$, $a \in A$.
So  $\frac{d}{dx}$ defines an $R$-derivation $R[x]\h\ra R[x]\h$ which sends $R[s]\h$ into itself. Furthermore $R[s]\h$ and $R[x]\h$ have the same constants with respect to this derivation (the constant ring is $R$ in both cases). Since $R[s]\h$ is integrally closed and $R[s]\h\subset R[x]\h$ is integral, and since $\frac{dx}{dx}=1\in R[s]\h$ it follows, by Lemma \ref{l2}  that $x\in R[s]\h$, hence $A=R[x]\h$, in other words the restriction homomorphism $\rho:\cO(J^1(E))\ra \cO(\widehat{{\mathbb G}_a})$ is an isomorphism.

\medskip

{\it Step 6}. Exactly as in Step 1 we get that the restriction map
\begin{equation}
\label{suna}
\cO(J^1(E)\times J^1(E))\ra \cO(\widehat{{\mathbb G}_a}\times \widehat{{\mathbb G}_a})=
\cO(\widehat{{\mathbb G}_a})\widehat{\otimes} \cO(\widehat{{\mathbb G}_a})\end{equation}
is injective. Since, by Step 5, the restriction homomorphism $\rho:\cO(J^1(E))\ra \cO(\widehat{{\mathbb G}_a})$ is an isomorphism we conclude that \ref{suna} is an isomorphism.
Let $\psi_1\in \cO(J^1(E))$ be the unique function whose restriction to $\widehat{{\mathbb G}_a}$ is $x$. Then $\psi_1$ defines a homomorphism $J^1(E)\ra \widehat{{\mathbb G}_a}$ (hence  a non-zero $\d$-character of order $1$ in the sense of \cite{book}, p. 190). By \cite{book}, p. 194, Proposition 7.15, $E$ has a lift of Frobenius, a contradiction. This ends the proof of assertion 1 in our Theorem.
\qed

\bigskip

{\it Proof of assertion 2}.

\medskip

{\it Step 1}. We start by assuming that $J^1(E)$ has a lift of Frobenius and we shall derive a contradiction. By the universality property of $p$-jets the lift of Frobenius on $J^1(E)$ induces a section $s:J^1(E)\ra J^1(J^1(E))$, in the category of $p$-formal schemes, of the canonical projection $\pi^1_{J^1(E)}:J^1(J^1(E))\ra J^1(E)$. We may (and will) assume $s(0)=0$.
Consider the homomorhism
$$l^{11}\d:J^1(J^1(E))\ra L_{\d}(E)=\widehat{{\mathbb G}_a}$$
defined in \ref{alpha}. By Corollary \ref{cartesian2} we have that $Ker(l^{11}\d)=J^2(E)$.
Now the composition 
$$l^{11}\d \circ s:J^1(E)\ra \widehat{{\mathbb G}_a}$$
defines a function in $\cO(J^1(E))$ sending $0$ into $0$. By the first assertion of Theorem \ref{main}, $l^{11}\d\circ s$ is a constant function, hence $s=0$, hence $s$ factors through $Ker(l^{11}\d)=J^2(E)$; in other words $s$ induces a section 
$\sigma:J^1(E)\ra J^2(E)$
of the canonical projection $\pi^2_E:J^2(E)\ra J^1(E)$; this is a section in the category of $p$-formal schemes (not a priori a group homomorphism). Also $\sigma(0)=0$.

\medskip

{\it Step 2}. Recall from \cite{char},  Proposition 2.2 and Lemma 2.3, that $Ker(\pi^2_E)$ is isomorphic to  $\widehat{{\mathbb G}_a}$.
Fix a point $P\in J^1(E)(R)$. The morphism of formal schemes
$f_P:J^1(E)\ra J^2(E)$ defined by 
$$f_P(Q)=\sigma(P+Q)-\sigma(P)-\sigma(Q)$$
factors through $Ker(\pi^2_E)$ hence gives rise to a morphism (still denoted by) 
$$f_P:J^1(E)\ra \widehat{{\mathbb G}_a}=\widehat{{\mathbb A}^1}.$$
 Since, by assertion 1) of the Theorem,  $\cO(J^1(E))=R$ we must have $f_P\in R$. Since $f_P(0)=0$ we get $f_P=0$. Since $P$ was arbitrary it follows  that $\sigma$ is a group homomorphism. Then the morphism
 $\chi:J^2(E)\ra J^2(E)$ defined by
 $$\chi(P)=P-\sigma(\pi^2_E(P))$$
 is a homomorphism. Clearly $\chi$ factors through a homomorphism (still denoted by) $\chi:J^2(E)\ra Ker(\pi^2_E)=\widehat{{\mathbb G}_a}$ and $Ker(\chi)=Im(\sigma)$ as formal schemes. In particular $\chi$ gives rise to a $\d$-character of order $2$ of $E$. So $\chi=c\cdot \psi$ for some $c\in R$, where $\psi=\psi_{2,can}$ is a canonical $\d$-character of order $2$ of $E$; cf. 
 \cite{book}, p. 201, Definition 7.24, and p. 205, Theorem 7.34.
 
 \medskip
 
 {\it Step 3}. 
 Let $Z$ be the subscheme of $J^2(Y)$ defined by the ideal generated by $\chi\in \cO(J^2(E))$. Then\begin{equation}
\label{yyy}
\overline{Z}\ra  \overline{J^1(E)},
\end{equation}
is an isomorphism. 
Let  $Y$ be an affine open subset of $E$ such that $\cO(Y)$ has an \'{e}tale coordinate
$T\in \cO(Y)$.
Now by 
\cite{book}, Theorem 7.22 plus Equation 7.73,  we have  identifications 
$$\cO(J^1(Y))=\cO(Y)[T']\h,\ \ \ \cO(J^2(Y))=\cO(Y)[T',T'']\h,$$ where $T',T''$ are new variables, $T'=\d T$, $T''=\d^2 T$, and 
$$\psi\in \cO(Y)[T']\h+p\cO(Y)[T',T'']\h.$$
 Hence  the reduction mod $p$, $\overline{\chi}$, of $\chi=c \psi$  belongs to $\cO(\overline{Y})[T']$ and so the map 
$$\cO(\overline{Y})[T']\ra \frac{\cO(\overline{Y})[T',T'']}{(\overline{\chi})}\simeq  \frac{\cO(\overline{Y})[T']}{(\overline{\chi})}[T'']$$
induced by \ref{yyy} is not an  isomorphism. This is a contradiction, which concludes the proof.
\qed

\end{document}